\documentclass[10pt]{article}
\usepackage[T1]{fontenc}              % norsk tegnsett (זרו)
\usepackage[latin1]{inputenc}  

\usepackage{amsmath,amsfonts,amssymb} 
\usepackage{amsthm}
\usepackage{graphicx}                 
\usepackage[round]{natbib}
\usepackage{geometry}
\usepackage{fullpage}
\usepackage{setspace}
\usepackage{epstopdf}
\usepackage{lineno}
\newcommand{\vc}[1]{\mathbf{#1}}
\newcommand{\mat}[1]{\mathbf{#1}}

\DeclareMathOperator{\toD}{\stackrel{d}{\to}}
\DeclareMathOperator{\toP}{\stackrel{P}{\to}}
\DeclareMathOperator{\E}{E}
\DeclareMathOperator{\diag}{diag}
\DeclareMathOperator{\var}{var}

\DeclareMathOperator{\cor}{cor}

\newtheorem{theorem}{Theorem}
\newtheorem{lemma}{Lemma}

\newtheorem{definition}{Definition}

\theoremstyle{remark}

\title{When and why are principal component scores a good tool for visualizing high-dimensional data?\\[1cm] \large Running headline: PC scores and high-dimensional data}

\author{Kristoffer Hellton, \\Department of Mathematics, University of Oslo\\Oslo Centre for Biostatistics and Epidemiology, Department of Biostatistics, University of Oslo \\\\Magne Thoresen \\
Oslo Centre for Biostatistics and Epidemiology, Department of Biostatistics, University of Oslo}

\begin{document}
\maketitle

\onehalfspacing 
\newpage
\begin{abstract}
Principal component analysis (PCA) is a popular dimension reduction technique often used to visualize high-dimensional data structures. In genomics, this can involve millions of variables, but only tens to hundreds of observations. Theoretically, such extreme high-dimensionality will cause biased or inconsistent eigenvector estimates, but in practice the principal component scores are used for visualization with great success. In this paper, we explore when and why the classical principal component scores can be used to visualize structures in high-dimensional data, even when there are few observations compared to the number of variables. Our argument is two-fold: First, we argue that eigenvectors related to pervasive signals will have eigenvalues scaling linearly with the number of variables. Second, we prove that for linearly increasing eigenvalues, the sample component scores will be scaled and rotated versions of the population scores, asymptotically. Thus the visual information of the sample scores will be unchanged, even though the sample eigenvectors are biased. In the case of pervasive signals, the principal component scores can be used to visualize the population structures, even in extreme high-dimensional situations. 
\end{abstract}


\section{Introduction}
Principal component analysis (PCA) is the workhorse of variable reduction in applied data analysis. It is used to construct a small number of informative scores from the original data, and these scores are then used for visualization of data or in conventional classification, clustering or regression methods. PCA is highly useful in the context of modern high-dimensional data analysis, where the number of measured variables, $p$, exceeds the sample size, $n$. Genomics is an application where data exploration is typically done by visually investigating the first few principal component (PC) scores. 

The asymptotic behavior of PCA in the high-dimensional setting has attracted a substantial amount of attention the last few years. It was shown by \citet{paul2007asymptotics} and \citet{johnstone2009consistency} that if $p >n$, the population eigenvalues and eigenvectors in PCA are not consistently estimated by the sample eigenvalues and eigenvectors. This was shown in an asymptotic framework where $p/n\to \gamma>0$ as both $p, n\to \infty$, while the population eigenvalues remained fixed. \citet{jung2009pca} introduced the High Dimension, Low Sample Size (HDLSS) asymptotic framework, where instead $n$ is fixed and the population eigenvalues scale with $p$ as $p\to\infty$. 

As a consequence of the eigenvector inconsistency, \citet{johnstone2009consistency} argued that the principal components should have a sparse representation, motivating the development of sparse PCA \citep{zou2006sparse,witten2009penalized}. \citet{cai2012minimax} stated: ``It is now well understood that in the high-dimensional setting the standard sample covariance matrix does not provide satisfactory performance and regularization is needed.'' 

In many genomic applications, however, classical (non-sparse) PCA and PC score plots are successfully used to discover and visualize population structures, suggesting that the scores are suitable for analysis in some situations. As \citet{lee2010convergence} noted: ``Inconsistency of the sample eigenvectors does not necessarily imply poor performance of principal component analysis''. We explain this paradoxical situation through the asymptotic behavior of the PC scores when $n$ is small and $p \to \infty$. Insights into the situation were given by \citet{shen2012high}, who showed that for a special case of the HDLSS framework the sample scores are scaled versions of the population scores, and by \citet{yang2014principal}, who proved that PC scores will correctly identify different strata in a stratified population. 

We first review the use of PCA with high-dimensional data and earlier asymptotic results. Next, we introduce the concept of pervasive effects and the connection to PCA, demonstrating how pervasive eigenvectors lead to population eigenvalues scaling linearly with the data dimension.  Two applications where pervasive structures are reasonable will be discussed. Further, under the assumption that the eigenvalues scale linearly with the data dimension, we derive the joint asymptotic limiting distribution of the sample principal component scores and show how the asymptotic distribution will impact the use of the component scores. 
In the case of two-dimensional visualizations, we show that pairs of sample component scores will be scaled and approximately rotated versions of the population scores with an additional error term. Based on the Taylor expansion of the distribution of the sample eigenvectors, the error term is seen to be negligible for large to moderate sample sizes, depending on the strength of the different signals. Thus the relative positions of the scores and thereby the visual information, will be preserved. 

\section{Principal component analysis}
\subsection{Methods and notation}
PCA reduces the dimension of a data matrix by constructing orthogonal linear combinations of the variables, such that the combinations express the directions with maximal variance in the data. The first principal component is the normalized linear combination of variables with the highest variance, while the second principal component will be the linear combination with the highest variance orthogonal to the first component, and so on. The mathematical basis of PCA is the eigendecomposition of the covariance matrix. 

Let  $\mat{X}= [\vc{x}_{1},\dots,\vc{x}_{n}]$ be a $p \times n$ data matrix, where $\vc{x}_{i}=[x_{i1}, \dots, x_{ip}]^{T}$ are independent and identically distributed with $\E \vc{x}_{i} = \vc{0}$ and $\var \vc{x}_{i}= \Sigma$. The eigendecomposition of the covariance matrix is given by
\[\Sigma = \mat{V} \mat{\Lambda} \mat{V}^{T}, \]
where $\mat{\Lambda}$ is the diagonal matrix of the eigenvalues $\lambda_{1} \geq \dots \geq \lambda_{p}$ and $\mat{V} = [\vc{v}_{1}, \dots, \vc{v}_{p}]$ is the matrix of eigenvectors. The population principal components are defined to be the linear combinations given by the eigenvectors of $\Sigma$, and the vector of the resulting scores is denoted
\begin{equation*}
\vc{s}_{j}^{T} = \vc{v}_{j}^{T}\mat{X} = [\vc{v}_{j}^{T}\vc{x}_{1}, \dots, \vc{v}_{j}^{T}\vc{x}_{n}], \quad j=1, \dots,p.
\end{equation*}
The variance of the component scores $\vc{s}_{j}^{T}$ is given by the $j$th eigenvalue, such that the vector of standardized population component scores is given by
\[\vc{z}_{j}^T = \frac{\vc{v}_{j}^{T}\mat{X}}{\sqrt{\lambda_{j}}},\quad j=1, \dots,p, \]
where the $j$th vector of scores is $\vc{z}_{j}^T=[z_{j1}, \dots, z_{jn}]$ and $\mat{Z} = [\vc{z}_{1},\dots,\vc{z}_{p}]$, an $n \times p$ matrix. 

In applied data analysis, the eigendecomposition is based on the sample covariance matrix with the proper centering, denoted by $\hat{\Sigma} = \frac{1}{n}\mat{X}\mat{X}^{T}$: 
\[\hat{\Sigma} =\hat{\mat{V}}\hat{\mat{\Lambda}}\hat{\mat{V}}^{T}. \]
Here $\hat{\mat{\Lambda}}=\diag(\hat{\lambda}_{1}, \dots, \hat{\lambda}_{p})$ contains the sample eigenvalues and $\hat{\mat{V}} = [\hat{\vc{v}}_{1},\dots,\hat{\vc{v}}_{p}]$ the corresponding sample eigenvectors. Following the earlier notation, the sample component scores are constructed as
\[\hat{\vc{s}}_{j}^{T}= \hat{\vc{v}}_{j}^{T}\mat{X}, \quad j=1, \dots,p,\]
and the sample standardized scores as
\[\hat{\vc{z}}_{j}^{T}= \frac{\hat{\vc{v}}_{j}^{T}\mat{X}}{\sqrt{\hat{\lambda}_{j}}}, \quad j=1, \dots,p.\]
For the rest of the paper, we assume that the population eigenvalues follow the spiked eigenvalue model introduced by \citet{johnstone2001distribution}, where the first $m$ population eigenvalues are substantially larger than the remaining non-spiked eigenvalues.

\subsection{Summary of earlier results}
\citet{anderson1963asymptotic} proved that the sample eigenvectors and -values,  $\hat{\vc{v}}$ and $\hat \lambda$, consistently estimate the population eigenvectors and -values, $\vc{v}$ and $\lambda$, when $p$ is fixed and $n \to \infty$. This result does not hold in the high-dimensional setting $(p \gg n)$. Starting with \citet{paul2007asymptotics} and \citet{johnstone2009consistency}, it has been shown that the sample eigenvalues and -vectors are not asymptotically consistent in the joint limit $p,n\to \infty$ when $p/n = \gamma >0$. \citet{paul2007asymptotics} proved that the inner product between the sample and the population eigenvector converges, when $\lambda_{j}> 1+\sqrt{\gamma}$, to
\[|\langle\hat{\vc{v}}_{j},\vc{v}_{j}\rangle|\to\sqrt{\left(1-\frac{\gamma}{(\lambda_{j}-1)^{2}}\right)/\left(1+\frac{\gamma}{\lambda_{j}-1}\right)}, \quad j=1,\dots,m.\]
%where the first $m$ eigenvalues are assumed to be substantially larger than the rest, as in the spiked covariance model. % introduced by \citet{johnstone2001distribution}. 
A different asymptotic framework began with the geometrical structure of high-dimensional data \citep{ahn2007high,hall2005geometric}. Building on this, \citet{jung2009pca} introduced the HDLSS regime, where $n$ is fixed and the $m$ spiked eigenvalues grow with the dimension $p$, according to 
$$\lambda_{j} = \sigma_{j}^{2}p^{\alpha}, \quad j =1, \dots,m.$$ 
In this asymptotic setting, consistency is defined in terms of $p$ asymptotics ($p \to \infty$) and the consistency of PCA depends on $\alpha$. Eigenvectors are estimated consistently when $\alpha > 1$, while the estimates are strongly inconsistent when $\alpha < 1$. In the boundary case $\alpha =1$, a situation explored by \citet{jung2012boundary}, the sample eigenvectors are neither consistent nor strongly inconsistent, but reach a limiting distribution depending on $n$. In the case of a single spiked eigenvalue,  $\lambda_1= \sigma_{1}^{2}p^{\alpha}$ and $\lambda_2=\cdots= \lambda_p= \tau^2$, and standard normally distributed population scores, the distribution of the sample eigenvector depends on $\alpha$ in the following way: 
\[|\langle\hat{\vc{v}}_{1},\vc{v}_{1}\rangle|\toD \left\{\begin{array}{ll}
1 &  \alpha > 1,\\
\left(1+\frac{\tau^{2}}{\sigma_1^{2}\chi^{2}_{n}}\right)^{-1/2} & \alpha = 1,\\
0 &  \alpha < 1, 
\end{array} \right.\]
where $\chi_{n}^{2}$ is a chi-square distribution with $n$ degrees of freedom. In the case of $\alpha = 1$, \citet{jung2012boundary} call $\hat{\vc{v}}_{1}$ inconsistent when $n$ is finite, but consistent when $n\to\infty$, and we will use the same terminology.

The main focus of the above-mentioned papers has been on the eigenvector inconsistency, and few results are concerned with principal component scores. One exception is \citet{lee2010convergence}, who established the asymptotic limit of the inner product between the sample and the population scores in the finite $\gamma$ regime of \citet{paul2007asymptotics}. They extended the results to prediction and found a theoretical asymptotic shrinkage factor for predicted scores. This shrinkage factor can be applied as a bias adjustment, which turns out to be useful in the context of genetic population stratification problems.  \citet{yata2009pca,yata2012effective} gave conditions for the sample scores to be consistent under the HDLSS framework in terms of the relation between $n$ and $p$ as $p, n\to \infty$. 

Further, \citet{shen2012high,shen2013surprising} investigated the ratio between the individual sample and population scores, instead of the inner product between the score vectors. They showed in the case of $\alpha > 1$ that the ratio $\hat{z}_{ij}/z_{ij}$ for $ i=1, \dots, n$, converges to a random  variable independent of $i$. This implies that asymptotically, a two dimensional plot of the $j$th  and $k$th PC scores is a scaled version of the population score plot. The visual information presented by the sample PC scores will therefore be the same as presented by the population scores. In this paper we investigate the same problem as \citet{shen2012high}, but in the situation where $\alpha=1$. We first motivate why this is an interesting assumption to make, and then derive the asymptotic behavior of the sample scores under this assumption.   

In summary, two different asymptotic frameworks have been utilized to evaluate the behavior of the sample eigenvalues and -vectors; either with fixed population eigenvalues, or with fixed sample size $n$. To understand when and why the classical PC scores are appropriate to use in a practical $p\gg n$ setting, we believe that the latter framework is more useful as the results from the HDLSS regime can be transferred to practical settings where $n$ is substantially smaller than $p$. The remaining challenge is to decide on an appropriate $\alpha$ for a specific application or real study. We will therefore show how the growth rate of the eigenvalues can be interpreted in terms of the data structure or data-generating mechanism, which is often easier to relate to. In our opinion, if one can argue that a given data-generating mechanism and corresponding growth rate is reasonable for a certain application, one can also justify relying on the asymptotic results from the HDLSS framework. %Our aim is to translate the assumption about the eigenvalues into an assumption regarding the latent structure and the data-generating mechanism, as this is . 

More recently, \citet{fan2015asymptotics} presented an asymptotic framework merging the finite $\gamma$ and the HDLSS frameworks, giving the necessary relation between $n$, $p$ and the $\lambda_i$'s for the eigenvectors and -values to be consistent. This result will be further explored in the discussion.

\section{Pervasiveness and eigenvalues}\label{datastructure}
In this section we demonstrate how assuming a certain data-generating mechanism, e.g. pervasive eigenvectors, leads to eigenvalues scaling linearly with the dimension, corresponding to the case of $\alpha =1$ in the HDLSS framework \citep{jung2012boundary}. In the econometrics literature, the concept of pervasiveness is commonly used for latent factor estimation \citep{bai2003inferential,fan2011high,fan2013large,bai2012efficient}. A pervasive factor is here thought to be an unobserved variable affecting most or a substantial proportion of the observed variables. In a situation where the dimension increases, pervasiveness can be formulated in terms of the asymptotic proportion of non-zero variable coefficients: 
\begin{definition}[Pervasiveness]\label{pervasiveness}
A sequence of $p$-dimensional vectors $\vc{v} = [v_{1}, \dots, v_{p}]^{T}$ fulfills the pervasiveness assumption if the proportion of non-zero entries $r_{p} = \frac{1}{p}\sum_{j=1}^{p}I_{\{v_{j}^{2}>0\}}$ fulfills:
\[\lim_{p \to\infty } r_{p} > 0.\]
\end{definition}
The definition of pervasiveness stands in contrast to the definition of \emph{sparsity} where the proportion of non-zero entries must converge to zero as $p \to \infty$, for instance by fixing the number of non-zero eigenvector coefficients.

\subsection{Pervasiveness and PCA}
The definition of pervasiveness is also relevant for PCA. In the simplest situation, the spiked covariance model with a single component \citep{johnstone2001distribution, nadler2008finite, johnstone2009consistency}, each observation is given by % a score $z_{i}$ and noise $\vc{\varepsilon}_{i}$ with zero mean and variance $\var (z_{i})=1$ and , 
\begin{equation*}
\vc{x}_{i} = \vc{v}z_{i} + \vc{u}_{i}, \qquad \var (\vc{u}_{i})=\sigma^{2}I_{p},\quad\var (z_{i})=1, \qquad  i=1, \dots, n.
\label{factorModel}
\end{equation*}
Then the population covariance matrix of $\vc{x}_{i}$ is given by 
\[\Sigma = \vc{v}\vc{v}^{T} + \sigma^{2} I_{p}.\]
The normalized vector $\vc{v}/ \|\vc{v}\|$ is an eigenvector of $\vc{v}\vc{v}^{T}$ and therefore also an eigenvector of $\Sigma$. The corresponding eigenvalue of $\vc{v}\vc{v}^{T}$ is given by the normalizing constant $\|\vc{v}\|= \sum_{j=1}^{p}v_{j}^{2}$. Thus the largest eigenvalue of $\Sigma$ is given by 
\[\lambda_{1} = \sum_{j=1}^{p}v_{j}^{2} + \sigma^{2}.\]
The dependence of the largest population eigenvalue of $\Sigma$ on the dimension, $p$, can therefore be determined by the coefficients of $\vc{v}$. 

If the values of $v_{1},\dots,v_{p}$ are constant (or independent of $p$) and fulfills the pervasiveness assumption, %as $p\to \infty$, 
there exist two constants $ c_{1} < c_{2}$, such that the largest population eigenvalue of the covariance matrix $\Sigma$ satisfies
\[c_{1}p + \sigma^{2} \leq \lambda_{1} \leq c_{2}p + \sigma^{2}. \]
The first eigenvalue will therefore scale linearly with the dimension as $p\to\infty$, while the other eigenvalues, $\lambda_{j}= \sigma^{2}$ for $j = 2, \dots, p$, will remain constant. 

The pervasive signal can also be interpreted in terms of the covariance matrix. Assume that $\Sigma$ is divided into the following independent equicorrelated blocks or sub-matrices:
\begin{equation*}
\Sigma= \begin{bmatrix}
\Sigma_{1} & 0 & \dots & 0 \\
0 &  \Sigma_{2}& 0 & \vdots\\
\vdots & 0&\Sigma_{3}  & 0 \\
0 & \dots &0 &  \sigma^{2}I \\
\end{bmatrix},  \quad \Sigma_{j} =\begin{bmatrix}
1+\sigma^{2} & \rho_{j}  & \dots & \rho_{j} \\
\rho_{j} & 1+\sigma^{2}  & &\vdots\\
\vdots &   &\ddots  & \rho_{j}\\
\rho_{j}&\hdots& \rho_{j} & 1+ \sigma^{2}
\end{bmatrix}, \quad j=1,2,3,
\label{}
\end{equation*} 
with $1>\rho_{1}> \rho_{2} > \rho_{3}>0$. Then each block will be represented by an eigenvector with equal, non-zero coefficients for the variables within the block and zeros for all others, and the covariance matrix decomposes as $$ \Sigma = \vc{v}_{1}\vc{v}_{1}^{T}+  \vc{v}_{2}\vc{v}_{2}^{T} + \vc{v}_{3}\vc{v}_{3}^{T} +\sigma^{2}I_{p}.$$ 
If the dimension of the block $\Sigma_{j}$ is $r_j p$, a non-zero proportion $r_{j}$ of the total number of variables as $p\to\infty$, and $1 \geq r_{1}+r_{2}+r_{3}> 0$, the covariance matrix $\Sigma$ will have three spiked eigenvalues
\[\lambda_{1} = \rho_{1}r_{1}p + 1+\sigma^{2}-\rho_{1}, \quad \lambda_{2} = \rho_{2}r_{2}p + 1+\sigma^{2}-\rho_{2},\quad \lambda_{3} = \rho_{3}r_{3} p + 1+\sigma^{2}-\rho_{3}, \]
all scaling linearly with the dimension $p$. % and the corresponding eigenvectors will be pervasive. 
The remaining  $p-3$ eigenvalues are constant and do not depend on $p$. Each spiked eigenvalue represents one block of variables and the importance of the block is determined by the proportion $r_{j}$ and the degree of correlation $\rho_{j}$ within the block. A pervasive eigenvector can therefore also be interpreted as a block of variables, where the dimension is a non-vanishing proportion of the total number of variables. 

\subsection{Applications}
The definition of a pervasive factor connects the data-generating mechanism behind the data to the asymptotic behavior of the eigenvalues. There are many applications where pervasive structures are reasonable, %suggesting that the assumption of eigenvalues scaling linearly with $p$, is reasonable, 
and we will present two examples. 

\subsubsection{Example 1: Single nucleotide polymorphisms}
Single nucleotide polymorphisms (SNPs) are genetic loci or markers having at least two alleles with an associated allelic frequency in a human population. The neutral theory of molecular evolution states that allele frequencies at most genetic loci change due to two stochastic processes; mutation and random drift. If the main variation in the data sample stems from differences between ethnic populations, random allelic drift is the main driver behind changes in the genetic markers. This will give many and randomly distributed differences, and when new markers are measured we expect a certain proportion to be informative with respect to ethnicity. This corresponds to our notion of a pervasive association between ethnicity and the SNPs, and as the association for each SNP is constant, the corresponding eigenvalue will scale linearly with number of included variables. Due to the random drift, the longer two populations have been separated, the larger proportion of differentiating SNPs one expects. This is seen in larger eigenvalues and stronger signals when comparing for instance European and Japanese populations, than when comparing subgroups of the Japanese population \citep{yamaguchi2008japanese}. 

\subsubsection{Example 2: Stock returns and macroeconomic measurements}
The concept of pervasive factors is, as previously mentioned, common in econometrics. \citet{harding2013estimating} demonstrates for instance the presence of pervasive factors in US stock returns using different methods. One can argue that important factors such as supply and demand shocks will impact a large number of firms simultaneously, such that a few underlying factors will have pervasive effects on the stock returns. In the case of macroeconomic measurements, \citet{stock2005implications} argue that the general state of the economy will be a pervasive factor impacting a range of economic and financial indicators. 

\section{Asymptotic results for sample scores}\label{asymptotic}
In this section we present asymptotic results for the behavior of the sample PC scores when the population eigenvalues scale linearly with the dimension. The asymptotic framework follows the HDLSS regime for the case where $\alpha = 1$, with the general conditions specified by \citet{jung2012boundary}:
\begin{itemize}
	\item[(C1)]  The standardized PC scores have finite fourth moments $\E(z_{ij}^4)<\infty$, and are uncorrelated, but possibly dependent, fulfilling the $\rho$-mixing condition. If the maximal correlation coefficient is defined as 
\[\rho(l) = \sup_{j,f,g}|\cor(f,g)|,\quad f \in L_{2}(\mathcal{F}_{-\infty}^j), g \in L_{2}(\mathcal{F}_{j+l}^\infty),\]
where $\mathcal{F}_{K}^L $ is the $\sigma$-field of events generated by the variables $\vc{z}_{j}, K \leq j \leq L$, the $\rho$-mixing condition is satisfied when $\rho(l) \to 0$, as $l\to \infty$. 
	\item[(C2)]	  The non-spiked eigenvalues  $\lambda_{m+1}, \dots, \lambda_{p}$ fulfill 
\[\frac{\sum_{i=m+1}^p \lambda_{i}^2}{\left(\sum_{i=m+1}^p \lambda_{i}\right)^2}\to 0,\quad \frac{1}{p}\sum_{i=m+1}^p \lambda_{i} \to \tau^2, \quad p \to \infty.\] 
\end{itemize}

Let the $n \times p$ matrix $\mat{Z}$ be the population standardized principal component scores. For an $m$ spiked model, the following notation for the scaled population scores is introduced: 
\[\tilde{\mat{Z}}_{m}= \left(\sigma_{1}\vc{z}_{1}, \dots, \sigma_{m}\vc{z}_{m}\right),\] 
%$\tilde{\mat{Z}}_{m}= \left(\sigma_{1}\vc{z}_{1}, \dots, \sigma_{m}\vc{z}_{m}\right)$, 
which gives the $m \times m$ matrix $\mat{W} = \frac{1}{n}\tilde{\mat{Z}}_{m}^{T}\tilde{\mat{Z}}_{m}$. The stochastic eigenvalues of the matrix are denoted by $\hat{\lambda}_{j}(\mat{W})$ and the stochastic eigenvectors are denoted by $\hat{\vc{v}}_{j}(\mat{W})$ for $j=1,\dots,m$. 

\begin{theorem}\label{scoreDistribution}
Under conditions $(C_{1})$, $(C_{2})$ and assuming the $m$ spiked eigenvalues to scale linearly with the dimension, %$\lambda_{1} = \sigma^{2}_{1}p,\dots,\lambda_{m}= \sigma_{m}^{2}p$,
\begin{equation}
\lambda_{1} = \sigma^{2}_{1}p,\quad\lambda_{2} = \sigma^{2}_{2}p, \quad \cdots, \quad \lambda_{m}= \sigma_{m}^{2}p,
\label{linearity}
\end{equation} 
where $m$ is fixed and $\sigma^{2}_{1}> \sigma^{2}_{2} > \cdots > \sigma_{m}^{2}$ are distinct, the $i$th vector of the $m$ first sample principal component scores converges in probability to 
\[
\renewcommand\arraystretch{1.5}
\begin{bmatrix} \hat{z}_{i1} \\  \vdots \\ \hat{z}_{im}
\end{bmatrix} \toP 
\renewcommand\arraystretch{1.5}
%\underbrace{
\begin{bmatrix}  
\sqrt{1/\hat{\lambda}_{1}(\mat{W})}  &&  0\\
 & \ddots & \\
0 &&  \sqrt{1/\hat{\lambda}_{m}(\mat{W})} 
\end{bmatrix}
%}_{Scaling}
%\underbrace{
\renewcommand\arraystretch{1.5}
\begin{bmatrix}  
&  & \\
\hat{\vc{v}}_{1}(\mat{W}) & \cdots & \hat{\vc{v}}_{m}(\mat{W}) \\
&  & \\
 \end{bmatrix}^T %}_{Rotation}
\begin{bmatrix} \sigma_{1} z_{i1}  \\ \vdots \\ \sigma_{m}z_{im}
\end{bmatrix}\]
where $\hat{\lambda}_{j}(\mat{W})$ and $\hat{\vc{v}}_{j}(\mat{W}) $ are the $j$th eigenvalues and eigenvectors of the stochastic matrix $\mat{W}= \frac{1}{n}\tilde{\mat{Z}}_{m}^{T}\tilde{\mat{Z}}_{m}$. The proof is given in the Appendix. 
% 
% \begin{align}
%   \hat{\mat{Z}}_{m} = \mat{V}_{m}(p^{-1}n\hat\Sigma_{D}) \toP \mat{V}\left(\tilde{\mat{Z}}_{m}\tilde{\mat{Z}}_{m}^T\right) = \hat{\mat{\Lambda}}(\tilde{\mat{Z}}_{m}^T\tilde{\mat{Z}}_{m})^{-1/2} \mat{V}(\tilde{\mat{Z}}_{m}^T\tilde{\mat{Z}}_{m})^T  \tilde{\mat{Z}}_{m},
% \end{align}
\end{theorem}

The first matrix, given by the sample eigenvalues, is a diagonal matrix and will act as a \emph{scaling matrix}. The second matrix, given by the $m$ sample eigenvectors, is an orthogonal matrix and will therefore act as an $m$-dimensional \emph{rotation matrix}. Thus, the $m$-dimensional sample score vector will be a scaled and rotated version of the population score vector in $m$-dimensional space. %However, 

When principal component scores are used to visualize high-dimensional data, it is common to display the $m$ first important sample scores in two-dimensional plots (or possibly in 3D).  These score plots can be used to compare observations, detect subgroups or identify outliers. Theorem \ref{scoreDistribution} shows that plots of the sample scores can still give valid information about the population scores, despite the inconsistent eigenvectors. First, we present the special cases of $m=1$ and $m=2$, and then the general result for $m>2$. 

For a single component, $\mat{W}$ will be a scalar such that the sample eigenvalue is $\hat{\lambda}_{1}(\mat{W})=\frac{\sigma_{1}^2}{n}\vc{z}_{1}^T\vc{z}_{1}$ and the eigenvector is simply 1. The first sample principal component score will therefore converge in probability to a scaled version of the population score 
\[\hat{z}_{i1} \toP  \sigma_{1}\sqrt{ 1 /\hat{\lambda}_{1}(\mat{W})}\;z_{i1} =\sqrt{\frac{n}{\vc{z}_{1}^T\vc{z}_{1}}}\;z_{i1}.\]
If the population scores are \emph{iid} normally distributed, $z_{ij}\sim \mathcal{N}(0,1)$, the first eigenvalue is distributed as $\hat{\lambda}_{1}(\mat{W})\sim\frac{\sigma_{1}^2}{n}\chi_{n}^2$, a scaled chi-squared distribution with $n$ degrees of freedom. The scaling of the sample score is then distributed as $R_{1} \sim \sqrt{n/\chi_{n}^2}$, the same distribution found by \citet{shen2012high} in the case where $\alpha >1$. 

For $m=2$, the sample scores converge in probability to
\[\renewcommand\arraystretch{1.5}
\begin{bmatrix} \hat{z}_{i1} \\ \hat{z}_{i2}
\end{bmatrix} \toP 
\underbrace{
\begin{bmatrix} \sqrt{1/\hat{\lambda}_{1}(\mat{W})} &0\\0& \sqrt{1/\hat{\lambda}_{2}(\mat{W})}  
\end{bmatrix}}_{\text{Scaling}}
\underbrace{
\begin{bmatrix}\hat v_{11}(\mat{W}) & \hat v_{12}(\mat{W})\\
\hat v_{21}(\mat{W}) & \hat v_{22}(\mat{W})
\end{bmatrix}}_{\text{Rotation}}
\begin{bmatrix} \sigma_1z_{i1} \\ \sigma_2z_{i2}\end{bmatrix}, \quad i=1, \dots, n.\]
As the eigenvector matrix is orthogonal and the eigenvalue matrix is diagonal, the sample score vector becomes a scaled and rotated version of the population score vector. Again, if the population scores are \emph{iid} normally distributed, $z_{ij}\sim \mathcal{N}(0,1)$, $\mat{W}$ will be a $2 \times 2$ Wishart distributed matrix
\[n\mat{W} \sim W_{2}\left(\begin{bmatrix}
\sigma_{1}^{2}& 0\\
0&\sigma_{2}^{2} \\
\end{bmatrix},n\right),\]
and the distribution of the scaling and rotation follows from this stochastic matrix.  

When there are more than two important components, one typically selects pairs of scores and plot these against each other. The distribution of a pair of the $j$th and $k$th sample principal component scores can for $m > 2$ be decomposed into two parts:
\begin{equation}
\renewcommand\arraystretch{1.5}
\begin{bmatrix} \hat{z}_{ij} \\ \hat{z}_{ik}
\end{bmatrix} \toP 
\renewcommand\arraystretch{1.5}
\underbrace{
\begin{bmatrix} \sqrt{1/\hat{\lambda}_{j}(\mat{W})} &0\\0& \sqrt{1/\hat{\lambda}_{k}(\mat{W})}  
\end{bmatrix}}_{\text{Scaling}}
\renewcommand\arraystretch{1.5}
\underbrace{
\begin{bmatrix}\hat v_{jj}(\mat{W}) &\hat v_{jk}(\mat{W})\\
\hat v_{kj}(\mat{W}) & \hat v_{kk}(\mat{W})
\end{bmatrix}}_{\text{Approximate rotation}}
\begin{bmatrix} \sigma_jz_{ij} \\ \sigma_kz_{ik}\end{bmatrix} + % \underbrace{
\begin{bmatrix} \varepsilon_{ij} \\ \varepsilon_{ik}
\end{bmatrix}%}_{Noise}
,\label{result}
\end{equation}
for $i=1, \dots, n$, where 
\begin{equation}
\varepsilon_{ij}=\sqrt{1/\hat{\lambda}_{j}(\mat{W})} \sum_{l=1,l\neq j,k}^m \sigma_{l}z_{il}\hat v_{jl}(\mat{W}), \quad \varepsilon_{ik}=\sqrt{1/\hat{\lambda}_{k}(\mat{W})} \sum_{l=1,l\neq j,k}^m \sigma_{l}z_{il}\hat v_{kl}(\mat{W}). \label{error} 
\end{equation}

\subsection{Distributions of the eigenvectors}
The key feature of Equation \eqref{result} is that the $\varepsilon_{ij}$- and $\varepsilon_{ik}$-term can be considered as noise compared to the scaling and approximate rotation matrix, for large enough $n$. However, for finite sample, the noise distribution depends on the complicated distribution of the eigenvectors. We therefore need to introduce the classical large sample ($n\to \infty$) distribution of the sample eigenvectors to simplify the derivation. This might seem contradictory as Theorem \ref{scoreDistribution} was derived under fixed $n$ and $p\to \infty$, but it is the most accessible way to characterize the noise distribution. The large sample results show how the parameters $\sigma_{1}^{2},\dots, \sigma_{m}^{2}$ and $m$ influence the distribution of $\varepsilon_{ij}$ and $\varepsilon_{ik}$.

The large sample distributions of eigenvectors and -values were established in the normal case by \citet{anderson1963asymptotic} and extended to the case of non-normal distributions by \citet{davis1977asymptotic}. Suppose the standardized population scores $z_{ij}$ for $ i=1, \dots, n,$ and $j=1, \dots, m,$ are \emph{iid} with $\E (z_{ij})=0$, $\var (z_{ij})= 1$ and finite fourth moment, characterized by the fourth order cumulants denoted $\kappa_{abcd} = \E(z_{ia}z_{ib}z_{ic}z_{id})-\E(z_{ia}z_{ib})\E(z_{ic}z_{id})-\E(z_{ia}z_{ic})\E(z_{ib}z_{id})-\E(z_{ia}z_{id})\E(z_{ib}z_{ic})$. Then the $m \times m$ matrix 
$$\mat{W} = \frac{1}{n}\tilde{\mat{Z}}_{m}^{T}\tilde{\mat{Z}}_{m}, \quad \tilde{\mat{Z}}_{m}= \left(\sigma_{1}\vc{z}_{1}, \dots, \sigma_{m}\vc{z}_{m}\right), $$ will have population eigenvalues $\sigma^{2}_{1}, \dots, \sigma^{2}_m$ and population eigenvectors, $\vc{e}_{j} =[0,\dots, 1,\dots,0]^{T}$, unit vectors with a single '1' at position $j$. By \citet{davis1977asymptotic} the sample eigenvectors of $\mat{W}$ converge in distribution to
\begin{equation*}
   n^{1/2}\left(\hat{\vc{v}}_{j}(\mat{W})- \vc{e}_{j}\right) \toD \mathcal{N}(0, \Psi), \text{ as } n\to\infty,
   \label{asympVector}
 \end{equation*}
where the entries of the asymptotic covariance matrix $ \Psi$ are
  \[\Psi_{l,l} = \frac{\kappa_{lljj}+\sigma_{l}^2\sigma_{j}^2}{(\sigma_{j}^2-\sigma_{l}^2)^2}, \quad \Psi_{l,l'} = \frac{\kappa_{ll'jj}}{(\sigma_{j}^2-\sigma_{l}^2)(\sigma_{j}^2-\sigma_{l'}^2)}, \quad l,l' =1, \dots, m; l,l' \neq j, \]
while the $j$th row and column of $\Psi$ only contain zeros; $\Psi_{l,j}=\Psi_{j,l}=0$ for all $l=1, \dots, m$. 

The sample eigenvalues converge in distribution to
\begin{equation*}
 n^{1/2}\left(\hat{\lambda}_{j}(\mat{W})-\sigma^{2}_{j}\right) \toD \mathcal{N}(0, 2\sigma^4_{j}+\kappa_{jjjj}), \text{ as } n\to\infty,
 \label{asympValue}
\end{equation*}
and $\hat{\vc{v}}_{j}(\mat{W})$ and $\hat{\lambda}_{j}(\mat{W})$ will be asymptotically independent.

In the case of normally distributed scores, $z_{ij} \sim N(0,1)$, all fourth order cumulants are zero and the joint distribution of the eigenvectors simplifies to the result of \citet{anderson1963asymptotic}:
 \[n^{1/2}\left(\hat{\vc{v}}_{j}(\mat{W})-\vc{e}_{j}\right) \toD \mathcal{N}(0, \Psi), \text{ as } n\to\infty,\] %are asymptotically normally distributed with mean 0 and covariance matrix
where 
 \[\Psi = \sum_{l=1, l\neq j}^m \frac{\sigma_{j}^2\sigma_{l}^2}{(\sigma_{j}^2-\sigma_{l}^2)^2}\; \vc{e}_{l}\vc{e}_{l}^T = \diag\left(\frac{\sigma_{j}^2\sigma_{1}^2}{(\sigma_{j}^2-\sigma_{1}^2)^2},\dots,0,\dots,\frac{\sigma_{j}^2\sigma_{m}^2}{(\sigma_{j}^2-\sigma_{m}^2)^2}\right).\]
The Delta method gives the asymptotic distribution of the eigenvalue and eigenvector part of Equation \eqref{error}:
\[ n^{1/2} \hat{\lambda}_{j}(\mat{W})^{-1/2}\hat v_{jl}(\mat{W}) \toD\mathcal{N} (0, \psi_{jl}),\]
where the asymptotic variance is given by the partial derivatives
\[\psi_{jl} = 2\sigma^4_{j}\left(\left.-\frac{1}{2}\hat{\lambda}_{j}(\mat{W})^{-\frac{3}{2}}\hat v_{jl}(\mat{W})\right|_{(\sigma_{j}^{2},0)}\right)^{2} +  \frac{\sigma_{j}^2\sigma_{l}^2}{(\sigma_{j}^2-\sigma_{l}^2)^2}\left(\left.\hat{\lambda}_{j}(\mat{W})^{-\frac{1}{2}}\right|_{(\sigma_{j}^{2},0)}\right)^{2}  = \frac{\sigma_{l}^2}{(\sigma_{j}^2-\sigma_{l}^2)^2}. 
\]
A single $z_{ij}$ and $\mat{W}$ are asymptotically independent, expressing the fact that a single observation has no influence when the sample size becomes infinite. Hence, as functions of $\mat{W}$, both $\hat{\vc{v}}_{j}(\mat{W})$ and $\hat{\lambda}_{j}(\mat{W})$ must be asymptotically independent of a single $z_{ij}$. The asymptotic variance of the noise term is then given as
\[\var (n^{1/2} \varepsilon_{ij}) = 
\sum_{l=1,l\neq j,k}^m \sigma_{l}^2 \var(z_{il}) \var (n^{1/2} \hat{\lambda}_{j}(\mat{W})^{-1/2}\hat v_{jl}(\mat{W})) = \sum_{l=1, l\neq j,k}^m\frac{\sigma_{l}^4}{(\sigma_{j}^2-\sigma_{l}^2)^2}, \quad n\to\infty,
\]
and the terms $\varepsilon_{ij}$ and $\varepsilon_{ik}$ are asymptotically independent. The variance of the noise terms will therefore depend on the signal strengths as

\begin{equation}
\var \left(\varepsilon_{ij}\right) =  \frac{1}{n}\sum_{l=1, l\neq j,k}^m \left( \frac{\sigma_{j}^2}{\sigma_{l}^{2}}-1\right)^{-2} + O\left(\frac{1}{n^2}\right).
\label{asympVar}
\end{equation}
Based on Equation \eqref{asympVar}, the main driver of the noise variation is the smallest difference $\sigma_{j}^2 - \sigma_{l}^{2}$, or eigengap, between the signal strengths of the components in question and all other components. The noise can therefore become large if one of the components in the pair has one or more neighboring components with similar signal strength, $\sigma_{l}^{2}$. In practice, one should therefore not plot components with similar eigenvalues separately, paired with other components. Instead one should try to plot the similar components together, for instance by extending the representation from two dimensions to three dimensions. The exact distribution of the noise for finite $n$, in the case of normally distributed data, is explored in Appendix \ref{simulation} using simulations. 

\section{Implications for visualization of sample scores}
When the noise part is negligible, Equation \eqref{result} shows that a two-dimensional plot of the estimated $j$th and $k$th sample principal component scores will be a scaled and approximately rotated version of the two-dimensional plot of the population principal component scores. Thus the relative positions of the population scores, and thereby the visual information, will be preserved by the sample scores. The behavior of the scaling matrix and the approximate rotation matrix will depend on the distribution of the random eigenvalues and -vectors, $\hat{\lambda}_{j}(\mat{W})$ and $\hat{\vc{v}}_{j}(\mat{W})$. Depending on the value of these quantities, the visualization of the sample scores compared to the population scores will behave differently. The behavior can be divided into three main cases: scaling, rotation or a so-called saddle point. 

\subsection{Case I: Scaling}
\begin{figure}%
\centering
\includegraphics[width=\columnwidth]{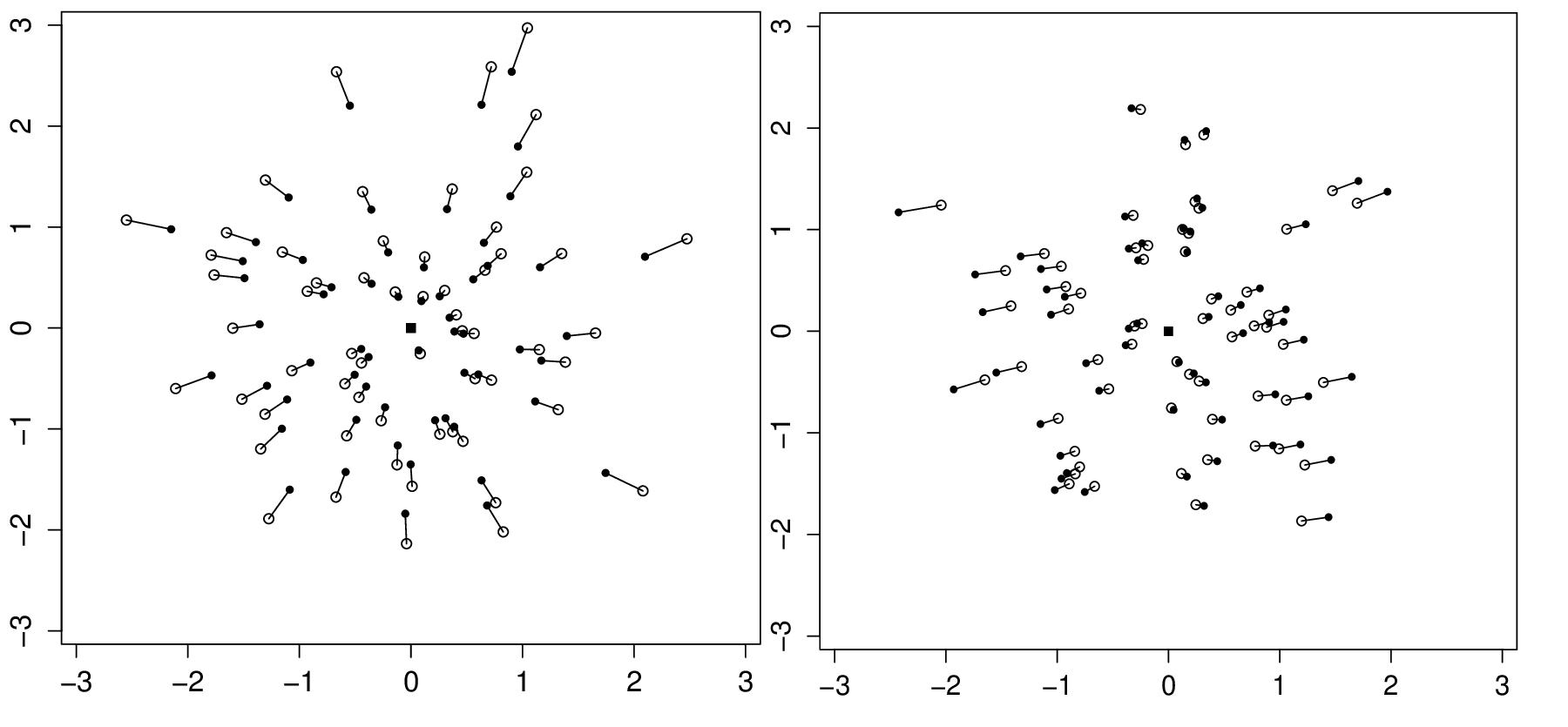}%
\caption{Two simulations of the first and second sample PC scores (dots) connected to the corresponding population PC scores (circles) with mostly radial scaling in positive and negative direction. The data are normally distributed with zero mean and identity covariance matrix with $p=12000$, $n=60$ and three large components ($m=3$): $\sigma^2_{1}=2.4, \sigma^2_{2}=0.8 $ and $ \sigma^2_{3} =0.1.$}

\label{RadialShift}%
\end{figure}

Firstly, if the scaling of both components is either significantly smaller or larger than 1, but there is no significant rotation, the sample scores will appear as a radial shift of the population scores. The shift can be either outwards or inwards from the origin as seen in Figure \ref{RadialShift}. This situation would be more common if the difference between the two signal strengths is large, such that the corresponding eigenvalues are well-separated. 

\begin{center}
	[Figure 1 here]
\end{center}

\subsection{Case II: Rotation}
\begin{figure}%
\centering
\includegraphics[width=\columnwidth]{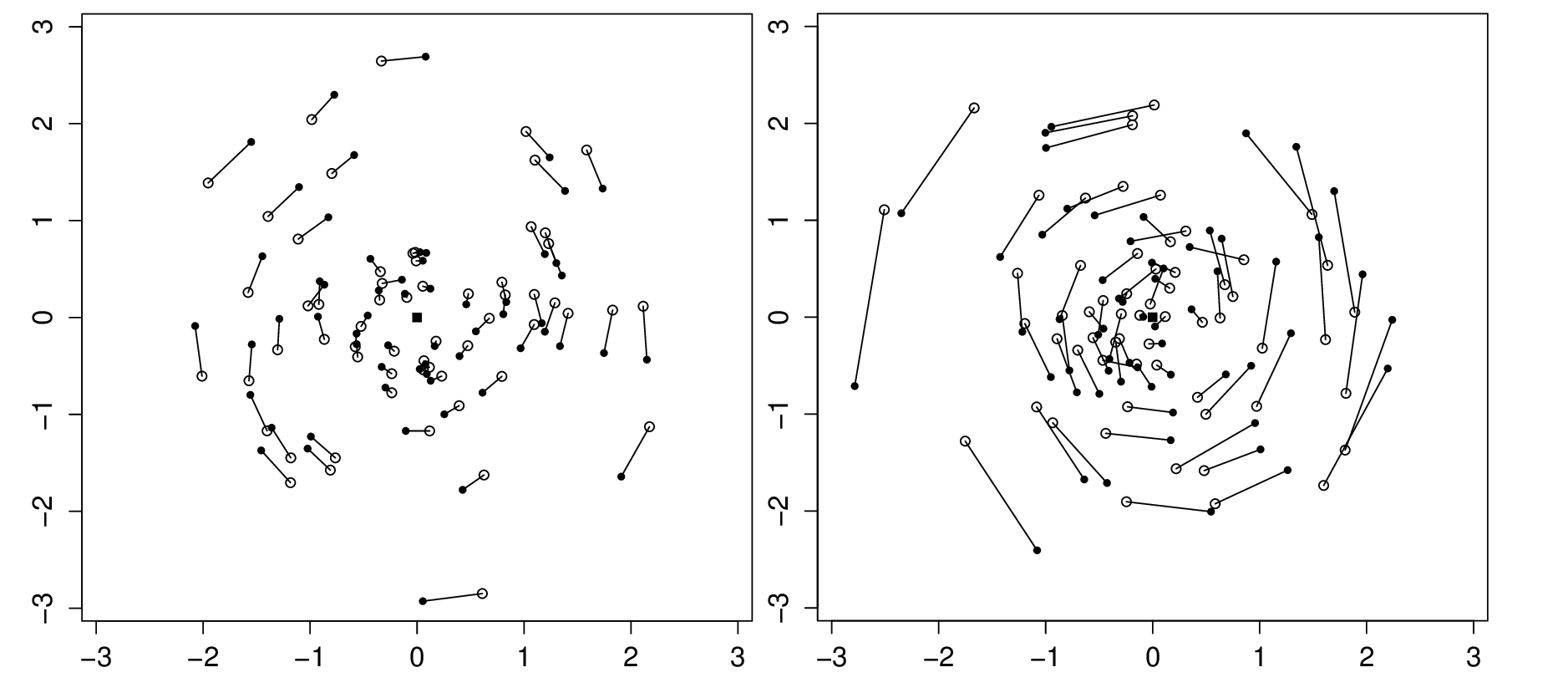}%
\caption{Two simulations of the first and second sample PC scores (dots) connected to the corresponding population PC scores (circles) with a large and small degree of rotation. The data are normally distributed with zero mean and identity covariance matrix with $p=12000$, $n=60$ and three large components ($m=3$): $\sigma^2_{1}=1.2, \sigma^2_{2}=0.8 $ and $ \sigma^2_{3} =0.4.$}%
\label{Rotation}%
\end{figure}

Secondly, if the approximate rotation is large, but the scaling in each direction is insignificant, the sample scores will appear as a rotation of the population scores. The rotation can be either in positive or negative angular direction around the origin as seen in Figure \ref{Rotation}. This situation would be more common if the two signal strengths are similar, such that the corresponding eigenvalues are not well-separated. Then, the estimated eigenvectors are mainly able to estimate the two-dimensional subspace spanned by the eigenvectors, and not the individual eigenvectors. 

\begin{center}
	[Figure 2 here]
\end{center}

\subsection{Case III: Saddle point}
\begin{figure}%
\centering
\includegraphics[width=\columnwidth]{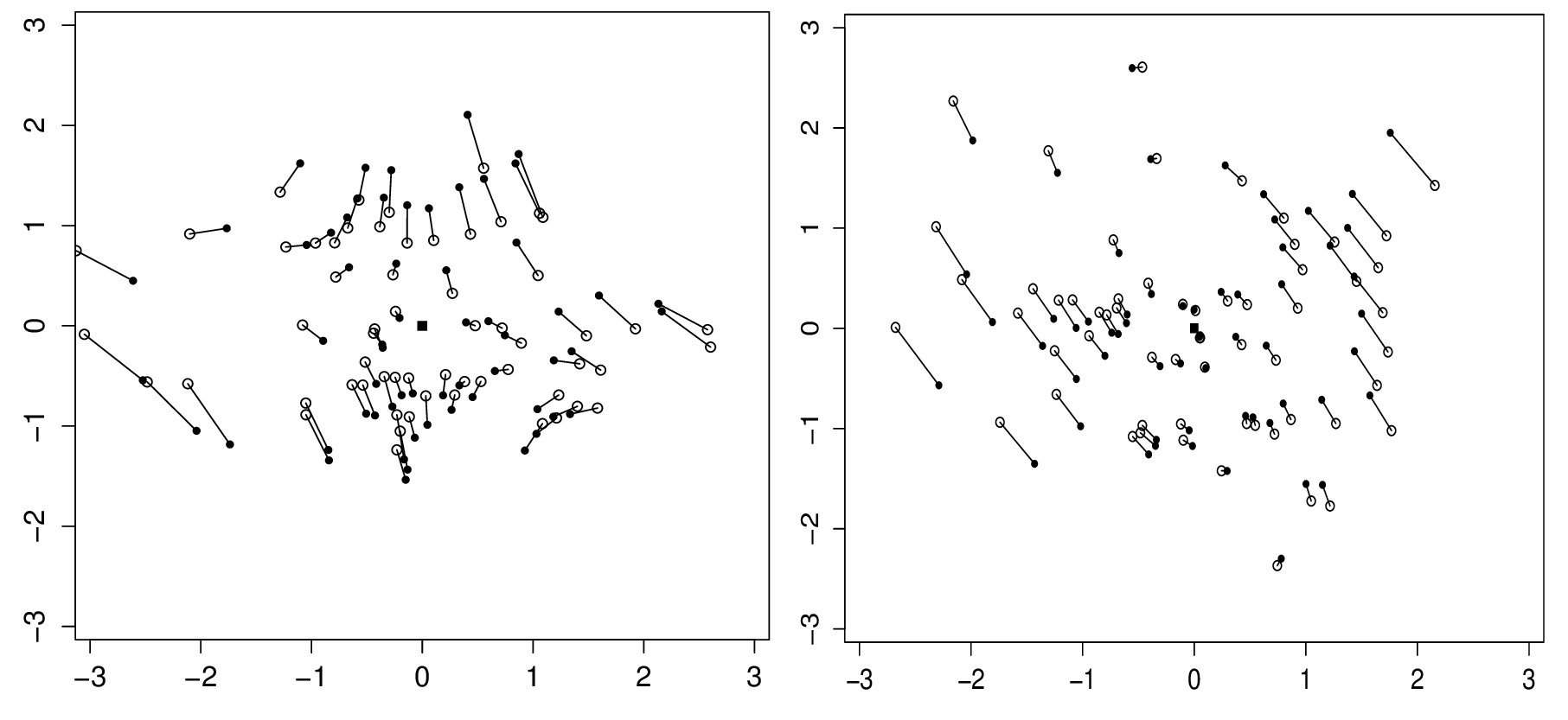}%
\caption{Two simulations of the first and second sample PC scores (dots) connected to the corresponding population PC scores (circles) where the scaling is less than 1 in one direction and larger than 1 in the other. The left panel shows a ``saddle point'' and the right shows an elongated rotation. The data are normally distributed with zero mean and identity covariance matrix with $p=12000$, $n=60$ and three large components ($m=3$): $\sigma^2_{1}=2.4, \sigma^2_{2}=0.8 $ and $ \sigma^2_{3} =0.1.$}%
\label{Saddle}%
\end{figure}

Lastly, we can have the special case seen in the first panel of Figure \ref{Saddle}, where the scaling is substantial, but smaller than 1 in one direction and larger than 1 in the other direction. An equivalent situation is found in the theory of linear dynamical systems, where it is termed a saddle point \citep[p. 66]{jordan2007nonlinear}. If the scaling is not large enough to overcome the rotation, one will instead see a fault-like rotation as in the second panel of Figure \ref{Saddle}.

\begin{center}
	[Figure 3 here]
\end{center}

\section{Discussion}
As \cite{paul2007asymptotics} and \cite{johnstone2009consistency} proved the eigenvectors to be inconsistent for fixed eigenvalues, the use of classical PCA for high-dimensional data became discredited, and this helped motivate the development of sparse PCA. With the contributions of \citet{shen2013surprising}, \citet{yang2014principal} and recently \citet{fan2015asymptotics}, however, the problems of classical PCA have become more nuanced. The important result of \cite{fan2015asymptotics} brings the different asymptotic frameworks together and establishes that the eigenvalues and -vectors will be consistent when $p/(n\lambda_j)\to 0$. Therefore, if  $n\to \infty$ and the population eigenvalues scale linearly with $p$ as in the case of pervasive signals, the eigenvectors, and hence PCA, will be consistent. This result can also be seen from Theorem \ref{scoreDistribution}, as the rotation and scaling matrix both converge to the identity matrix when $n\to \infty$. 

We show in addition that the visual information conveyed by the scores can be unbiased, or consistent in the terminology of \citet{jung2012boundary}, for finite $n$, and this will impact the use of PCA for practical data analysis. It shows that visual inspection of the PC scores will be a good way to assess population structures, even when $p \gg n$ as is commonly done in genomics. 

The remaining question is whether the pervasive signal model is reasonable for a certain application, and we have highlighted the example of ethnicity and SNP markers, an application where classical PCA is widely used. Further work also needs to explore how these results will affect clustering and regression, as these methods often use PC scores as inputs. 

\section*{Acknowledgments}
The authors would like to thank the editor and referees for their valuable comments, and J. S. Marron for fruitful discussions. The work was supported by the Norwegian Cancer Society. 

\bibliographystyle{sjs}
\bibliography{bibBiom}

\begin{thebibliography}{29}
\expandafter\ifx\csname natexlab\endcsname\relax\def\natexlab#1{#1}\fi

\bibitem[{Acker(1974)}]{acker1974absolute}
Acker, A.~F. (1974).
\newblock Absolute continuity of eigenvectors of time-varying operators.
\newblock \emph{Proceedings of the American Mathematical Society} \textbf{42},
  198--201.

\bibitem[{Ahn \emph{et~al.}(2007)Ahn, Marron, Muller \& Chi}]{ahn2007high}
Ahn, J., Marron, J.~S., Muller, K.~M. \& Chi, Y. (2007).
\newblock The high-dimension, low-sample-size geometric representation holds
  under mild conditions.
\newblock \emph{Biometrika} \textbf{94}, 760--766.

\bibitem[{Anderson(1963)}]{anderson1963asymptotic}
Anderson, T. (1963).
\newblock Asymptotic theory for principal component analysis.
\newblock \emph{Ann. Statist.} \textbf{34}, 122--148.

\bibitem[{Bai(2003)}]{bai2003inferential}
Bai, J. (2003).
\newblock Inferential theory for factor models of large dimensions.
\newblock \emph{Econometrica} \textbf{71}, 135--171.

\bibitem[{Bai \& Liao(2016)}]{bai2012efficient}
Bai, J. \& Liao, Y. (2016).
\newblock Efficient estimation of approximate factor models via penalized
  maximum likelihood.
\newblock \emph{Journal of Econometrics} \textbf{1}, 1--18.

\bibitem[{Cai \& Zhou(2012)}]{cai2012minimax}
Cai, T.~T. \& Zhou, H.~H. (2012).
\newblock Minimax estimation of large covariance matrices under l1 norm.
\newblock \emph{Statist. Sinica} \textbf{22}, 1319--1349.

\bibitem[{Davis(1977)}]{davis1977asymptotic}
Davis, A.~W. (1977).
\newblock Asymptotic theory for principal component analysis: Non-normal case.
\newblock \emph{Australian Journal of Statistics} \textbf{19}, 206--212.

\bibitem[{Fan \emph{et~al.}(2011)Fan, Liao \& Mincheva}]{fan2011high}
Fan, J., Liao, Y. \& Mincheva, M. (2011).
\newblock High dimensional covariance matrix estimation in approximate factor
  models.
\newblock \emph{Ann. Statist.} \textbf{39}, 3320.

\bibitem[{Fan \emph{et~al.}(2013)Fan, Liao \& Mincheva}]{fan2013large}
Fan, J., Liao, Y. \& Mincheva, M. (2013).
\newblock Large covariance estimation by thresholding principal orthogonal
  complements.
\newblock \emph{J. R. Statist. Soc. B} \textbf{75}, 603--680.

\bibitem[{Fan \& Wang(2015)}]{fan2015asymptotics}
Fan, J. \& Wang, W. (2015).
\newblock Asymptotics of empirical eigen-structure for ultra-high dimensional
  spiked covariance model.
\newblock Preprint arXiv:1502.04733.

\bibitem[{Hall \emph{et~al.}(2005)Hall, Marron \& Neeman}]{hall2005geometric}
Hall, P., Marron, J.~S. \& Neeman, A. (2005).
\newblock Geometric representation of high dimension, low sample size data.
\newblock \emph{J. R. Statist. Soc. B} \textbf{67}, 427--444.

\bibitem[{Harding(2013)}]{harding2013estimating}
Harding, M. (2013).
\newblock Estimating the number of factors in large dimensional factor models.
\newblock Unpublished manuscript.

\bibitem[{Johnstone(2001)}]{johnstone2001distribution}
Johnstone, I.~M. (2001).
\newblock On the distribution of the largest eigenvalue in principal components
  analysis.
\newblock \emph{Ann. Statist.} \textbf{29}, 295--327.

\bibitem[{Johnstone \& Lu(2009)}]{johnstone2009consistency}
Johnstone, I.~M. \& Lu, A.~Y. (2009).
\newblock On consistency and sparsity for principal components analysis in high
  dimensions.
\newblock \emph{J. Amer. Statist. Assoc.} \textbf{104}, 682--693.

\bibitem[{Jordan \& Smith(2007)}]{jordan2007nonlinear}
Jordan, D.~W. \& Smith, P. (2007).
\newblock \emph{Nonlinear ordinary differential equations: an introduction for
  scientists and engineers}.
\newblock New York.

\bibitem[{Jung \& Marron(2009)}]{jung2009pca}
Jung, S. \& Marron, J.~S. (2009).
\newblock {PCA} consistency in high dimension, low sample size context.
\newblock \emph{Ann. Statist.} \textbf{37}, 4104--4130.

\bibitem[{Jung \emph{et~al.}(2012)Jung, Sen \& Marron}]{jung2012boundary}
Jung, S., Sen, A. \& Marron, J.~S. (2012).
\newblock Boundary behavior in high dimension, low sample size asymptotics of
  {PCA}.
\newblock \emph{J. Mult. Anal.} \textbf{109}, 190--203.

\bibitem[{Lee \emph{et~al.}(2010)Lee, Zou \& Wright}]{lee2010convergence}
Lee, S., Zou, F. \& Wright, F.~A. (2010).
\newblock Convergence and prediction of principal component scores in
  high-dimensional settings.
\newblock \emph{Ann. Statist.} \textbf{38}, 3605.

\bibitem[{Nadler(2008)}]{nadler2008finite}
Nadler, B. (2008).
\newblock Finite sample approximation results for principal component analysis:
  A matrix perturbation approach.
\newblock \emph{Ann. Statist.} \textbf{36}, 2791--2817.

\bibitem[{Paul(2007)}]{paul2007asymptotics}
Paul, D. (2007).
\newblock Asymptotics of sample eigenstructure for a large dimensional spiked
  covariance model.
\newblock \emph{Statist. Sinica} \textbf{17}, 1617--1642.

\bibitem[{Shen \emph{et~al.}(2012)Shen, Shen, Zhu \& Marron}]{shen2012high}
Shen, D., Shen, H., Zhu, H. \& Marron, J.~S. (2012).
\newblock High dimensional principal component scores and data visualization.
\newblock Preprint arXiv:1211.2679.

\bibitem[{Shen \emph{et~al.}(2013)Shen, Shen, Zhu \&
  Marron}]{shen2013surprising}
Shen, D., Shen, H., Zhu, H. \& Marron, J.~S. (2013).
\newblock Surprising asymptotic conical structure in critical sample
  eigen-directions.
\newblock Preprint arXiv:1303.6171.

\bibitem[{Stock \& Watson(2005)}]{stock2005implications}
Stock, J.~H. \& Watson, M.~W. (2005).
\newblock Implications of dynamic factor models for {VAR} analysis.
\newblock Tech. rep., National Bureau of Economic Research.

\bibitem[{Witten \emph{et~al.}(2009)Witten, Tibshirani \&
  Hastie}]{witten2009penalized}
Witten, D.~M., Tibshirani, R. \& Hastie, T. (2009).
\newblock A penalized matrix decomposition, with applications to sparse
  principal components and canonical correlation analysis.
\newblock \emph{Biostatistics} \textbf{10}, 515--534.

\bibitem[{Yamaguchi-Kabata \emph{et~al.}(2008)Yamaguchi-Kabata, Nakazono,
  Takahashi, Saito, Hosono, Kubo, Nakamura \& Kamatani}]{yamaguchi2008japanese}
Yamaguchi-Kabata, Y., Nakazono, K., Takahashi, A., Saito, S., Hosono, N., Kubo,
  M., Nakamura, Y. \& Kamatani, N. (2008).
\newblock Japanese population structure, based on snp genotypes from 7003
  individuals compared to other ethnic groups: effects on population-based
  association studies.
\newblock \emph{The American Journal of Human Genetics} \textbf{83}, 445--456.

\bibitem[{Yang \emph{et~al.}(2014)Yang, Doksum \& Tsui}]{yang2014principal}
Yang, F., Doksum, K. \& Tsui, K.-W. (2014).
\newblock Principal component analysis {(PCA)} for high-dimensional data. {PCA}
  is dead. {L}ong live {PCA}.
\newblock \emph{Perspectives on {B}ig {D}ata Analysis: Methodologies and
  Applications} \textbf{622}, 1--22.

\bibitem[{Yata \& Aoshima(2009)}]{yata2009pca}
Yata, K. \& Aoshima, M. (2009).
\newblock {PCA} consistency for non-{G}aussian data in high dimension, low
  sample size context.
\newblock \emph{Communications in Statistics--Theory and Methods} \textbf{38},
  2634--2652.

\bibitem[{Yata \& Aoshima(2012)}]{yata2012effective}
Yata, K. \& Aoshima, M. (2012).
\newblock Effective {PCA} for high-dimension, low-sample-size data with noise
  reduction via geometric representations.
\newblock \emph{J. Mult. Anal.} \textbf{105}, 193--215.

\bibitem[{Zou \emph{et~al.}(2006)Zou, Hastie \& Tibshirani}]{zou2006sparse}
Zou, H., Hastie, T. \& Tibshirani, R. (2006).
\newblock Sparse principal component analysis.
\newblock \emph{J. Comp. Graph. Stat.} \textbf{15}, 265--286.

\end{thebibliography}

\vspace{\baselineskip}
\noindent Kristoffer Hellton, Dept. of Mathematics, University of Oslo, P.O.Box 1053 Blindern, NO-0316 Oslo, Norway. \\
E-mail: kristohh@math.uio.no

\appendix
\section*{Appendix}
\section{Proof of Theorem \ref{scoreDistribution}}
\begin{proof}
The proof follows from using the dual of the sample covariance matrix and the following lemma:
\begin{lemma}[\citet{jung2012boundary}]\label{lawLargeNumbers}
If assumptions (C1) and (C2) in Section \ref{asymptotic} hold, it follows that 
\[p^{-1} \sum_{j=m+1}^p\lambda_{j}\vc{z}_{j}\vc{z}_{j}^T \toP \tau^2 I_{n}, \quad \text{as } p \to \infty.\]
\end{lemma}
\noindent The dual of $\hat\Sigma$ is defined as the $n \times n$ symmetric matrix
$$  \hat\Sigma_{D} = \frac{1}{n} \mat{X}^T \mat{X}.$$
The dual has the singular value decomposition $ \hat\Sigma_{D} = \hat{\mat Z}_{n}\hat{\mat{\Lambda}}_{n}\hat{\mat Z}_{n}^T$ with $n$ eigenvalues $\hat{\mat{\Lambda}}_{n} =\diag(\hat{\lambda}_{1}, \dots, \hat{\lambda}_{n})$ and eigenvectors $\hat{\mat Z}_{n} = [\hat{\vc{z}}_{1},\dots,\hat{\vc{z}}_{n}],\hat{\vc{z}}_{j} = [\hat z_{j1}, \dots, \hat z_{jn}]^T $ corresponding to the $n$ first standardized sample principal components. With the $n \times m$ matrix
\[\tilde{\mat{Z}}_{m}= \left(\sigma_{1}\vc{z}_{1}, \dots, \sigma_{m}\vc{z}_{m}\right),\] 
we have the following two matrices: $\mat{W} = \frac{1}{n}\tilde{\mat{Z}}_{m}^T \tilde{\mat{Z}}_{m}$ and $\mat{W}_{D}= \frac{1}{n}\tilde{\mat{Z}}_{m} \tilde{\mat{Z}}_{m}^T$. 

By scaling the dual matrix $\hat\Sigma_{D}$ by $p^{-1}$ and inserting the $m$ first population eigenvalues from Equation \eqref{linearity}, one obtains the decomposition 
\begin{align*}
p^{-1}n \hat\Sigma_{D}&= \sum_{j=1}^m p^{-1}\lambda_{j} \vc{z}_{j}\vc{z}_{j}^T + p^{-1} \sum_{j=m+1}^p\lambda_{j}\vc{z}_{j}\vc{z}_{j}^T =  \sum_{j=1}^m \sigma_{j}^2 \vc{z}_{j}\vc{z}_{j}^T + p^{-1} \sum_{j=m+1}^p\lambda_{j}\vc{z}_{j}\vc{z}_{j}^T, \\
&=  n\mat{W}_{D} + p^{-1}\sum_{j=m+1}^p\lambda_{j}\vc{z}_{j}\vc{z}_{j}^T.
\end{align*}
Following Lemma 1 and scaling by $1/n$, the matrix converges in probability as $p\to \infty$ to 
\begin{equation}
 p^{-1} \hat\Sigma_{D} \toP \mat{W}_{D} + \frac{\tau^2}{n} I_{n},\label{converge}
\end{equation}
as $\mat{W}_{D}$ does not depend on the dimension $p$. The random variables $\mat{W}$ and $\mat{W}_{D}$ share the same non-zero eigenvalues
 $$\hat{\lambda}_{j}(\mat{W}) = \hat{\lambda}_{j}(\mat{W}_{D}), \text{ for } j = 1, \dots, m,$$
while the eigenvectors are proportional as
$$ \hat{\vc{v}}_{j}(\mat{W}_{D}) = \frac{\tilde{\mat{Z}}_{m}\hat{\vc{v}}_{j}(\mat{W})}{\sqrt{\hat{\lambda}_{j}(\mat{W})}}, \text{ for } j = 1, \dots, m,$$
or in matrix notation: $$\hat{\mat{V}}_{m}(\mat{W}_{D}) = \tilde{\mat{Z}}_{m} \hat{\mat{V}}(\mat{W}) \hat{\mat{\Lambda}}(\mat{W})^{-1/2}.$$ 
%The eigenvalues of a symmetric matrix $\mat{A}$ are continuous functions of the elements of $\mat{A}$ \citep{kato2013perturbation}. 
If the eigenvalues of a symmetric matrix $\mat{A}$ are distinct, the corresponding eigenvectors can be chosen such that they are continuous functions of the elements of $\mat{A}$ \citep{acker1974absolute}. From Equation \eqref{converge}, we then have that
$$ \hat{\vc{v}}_{j}(p^{-1}\hat\Sigma_{D}) \toP \hat{\vc{v}}_{j}\left(\mat{W}_{D} + \frac{\tau^2 }{n}I_{n}\right) = \hat{\vc{v}}_{j}\left(\mat{W}_{D}\right), \qquad j =1, \dots, m,$$
as the addition of an identity matrix does not change the eigenvectors. 
The sample principal components are given by the eigenvectors of the dual matrix, such that the result follows
\begin{align*}
  \hat{\mat{Z}}_{m} = \hat{\mat{V}}_{m}(p^{-1}\hat\Sigma_{D}) \toP \hat{\mat{V}}\left(\frac{1}{n}\tilde{\mat{Z}}_{m}\tilde{\mat{Z}}_{m}^T\right) =  \tilde{\mat{Z}}_{m} \hat{\mat{V}}(\mat{W}) \hat{\mat{\Lambda}}(\mat{W})^{-1/2}, 
  \end{align*}
or with the observations as columns:
\begin{align*}
  \hat{\mat{Z}}_{m}^T = \hat{\mat{\Lambda}}(\mat{W})^{-1/2} \hat{\mat{V}}^T(\mat{W})  \tilde{\mat{Z}}_{m}^T. 
\end{align*}
\end{proof}

\section{Simulation in the case of normally distributed data}\label{simulation}
To explore the distribution of $\varepsilon_{ij}$ in Equation \eqref{result} in a finite sample situation, we simulate normally distributed data and explore the standard deviation of the noise distribution. Particularly we investigate how the noise depends on the sample size, $n$, and signal strengths, $\sigma^{2}_{1}, \dots,\sigma_{m}^{2}$. % and the total number of spiked components $m$.

\begin{figure}%
\centering
\includegraphics[height=11cm,width=11cm]{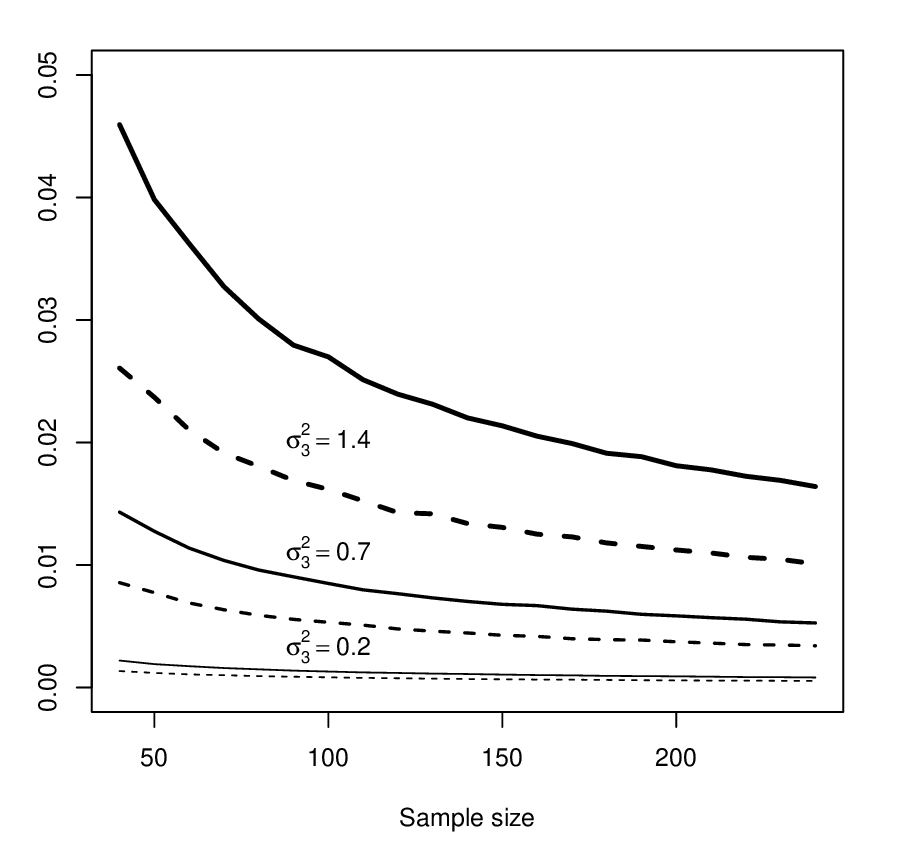}%
\caption{For $m=5$; the standard deviation of the noise $\varepsilon_{i1}$ (dashed curve) and $\varepsilon_{i2}$ (solid curve) as a function of sample size and three different ratios between the signal strengths of the second and third component (thick, medium and thin lines), based on 10 000 simulations of normally distributed variables.} \label{GraphSD}%
\end{figure}

We restrict the simulations to $m=5$ spiked components with the signal strengths
\[\sigma^2_1 = 12,\quad \sigma^2_2 = 8, \quad\sigma^2_4=0.1, \quad\sigma^2_5=0.02,\]
while $\sigma_{3}^{2}$ varies between $0.2$ and $1.4$, mainly evaluating the first and second component scores $j=1$ and $ k=2$. Figure \ref{GraphSD} shows the standard deviation of the noise distribution over 10000 simulations, demonstrating how the standard deviation decreases for increasing sample size with the asymptotic rate $n^{-1/2}$. The standard deviations of $\varepsilon_{i1}$ and $\varepsilon_{i2}$ are shown by dashed and solid lines, respectively, for different $\sigma_{3}^{2}$.

\begin{center}
	[Figure 4 here]
\end{center}

The main driver behind the variation is the smallest ratio between $\sigma_{1}^{2}$ and $\sigma_{2}^{2}$ and all the other signal strengths, the ratio $\sigma^2_{2}/\sigma^2_3$. This is seen when the standard deviation of $\varepsilon_{i1}$ and $\varepsilon_{i2}$ approximately doubles, when $\sigma^2_3$ changes from 0.7 to 1.4. As seen by Equation \eqref{asympVar}, the 4:3 ratio between $\sigma_{1}^{2}$ and $\sigma_{2}^{2}$ is also reflected in the standard deviation, as the ratio between the standard deviation of $\varepsilon_{i1}$ and $\varepsilon_{i2}$ is also approximately 4:3.  

\begin{figure}%
\centering
\includegraphics[height=11cm,width=11cm]{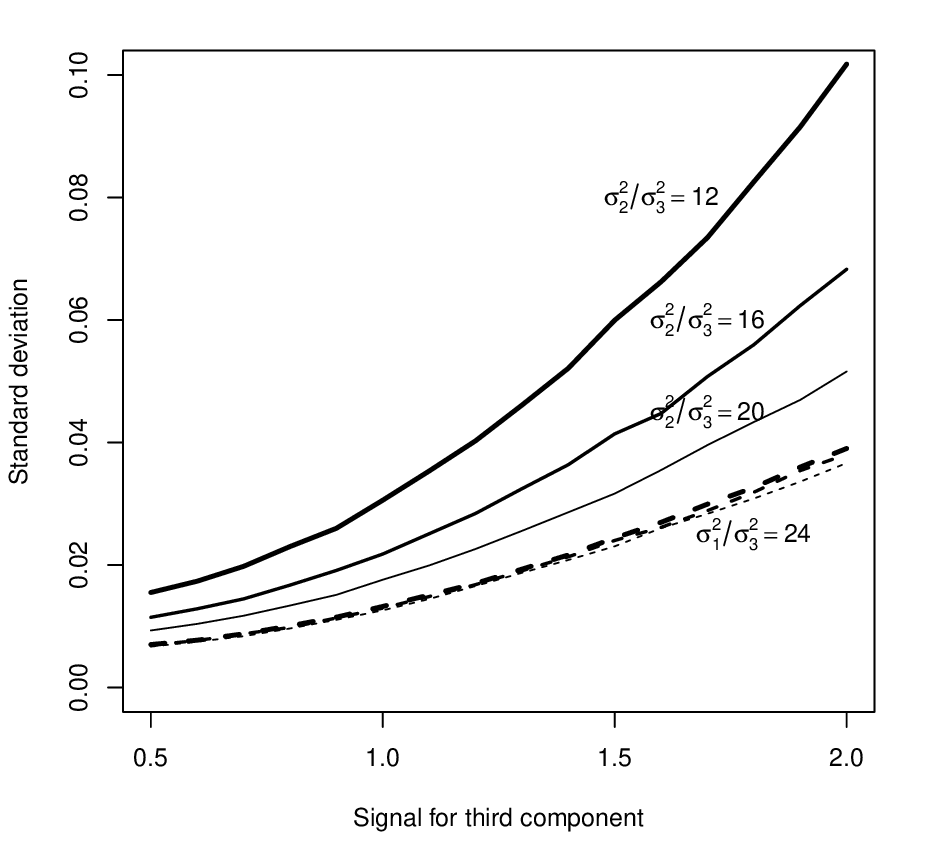}%
\caption{For $m=5$ and $n=60$; the standard deviation of the noise $\varepsilon_{i1}$ (dashed curve) and $\varepsilon_{i2}$ (solid curve) as a function of the signal strength of the third component, based on 10 000 simulations of normally distributed variables. The thick, medium and thin lines show three different ratios between the signal strengths of the second and third component.} \label{GraphSD2}%
\end{figure}

Figure \ref{GraphSD2} shows standard deviation of noise over 10000 simulations as a function of $\sigma_{3}^{2}$ and fixed sample size, $n=60$. The signal strength of the second component, $\sigma_{2}^{2}$, has three different values, shown by thick, medium and thin lines. The standard deviation of both $\varepsilon_{i1}$ (dashed line) and $\varepsilon_{i2}$ (solid line) will increase when $\sigma_{3}^{2}$ increases. When $\sigma_{2}^{2}$ has a smaller value, the standard deviation of $\varepsilon_{i2}$ increases more as $\sigma_{3}^{2}$ increases, and vice versa when $\sigma_{2}^{2}$ is larger. In both cases, the standard deviation of $\varepsilon_{i1}$ remains unchanged as long as the value of $\sigma_{1}^{2}$ is fixed.  

\begin{center}
	[Figure 5 here]
\end{center}

\clearpage

\end{document}